\documentclass{amsart}
\usepackage{amsmath, amssymb}
\usepackage{fullpage}
\usepackage{setspace}
\usepackage{hhline}

\newtheorem{theorem}{Theorem}

\newtheorem{proposition}[theorem]{Proposition}

\begin{document}

\title[Voting, the Symmetric Group, and Representation Theory]{Voting, the Symmetric Group, and Representation Theory}
\author{Zajj Daugherty, Alexander K. Eustis, Gregory Minton, and Michael E. Orrison}
\date{December 17, 2007}

\maketitle



\begin{abstract}
We show how voting may be viewed naturally from an algebraic perspective by viewing voting profiles as elements of certain well-studied $\mathbb{Q}S_n$-modules.  By using only a handful of simple combinatorial objects (e.g., tabloids) and some basic ideas from representation theory (e.g., Schur's Lemma), this allows us to recast and extend some well-known results in the field of voting theory.  
\end{abstract}


\section{Introduction}

For more than 25 years, Donald Saari has been systematically developing a powerful geometric approach to understanding, explaining, and constructing paradoxes that occur in voting \cite{merlin-saari-1997-copeland-2}-\cite{saari-valognes-1998-geometry}.  One of the keys to Saari's geometric approach is the fact that the collection of votes from an election can often be encoded naturally as a vector, which we call a \emph{profile}, and that election procedures can often be viewed as, or are related to, linear transformations.

By focusing on specific geometric structures in this vector space setting, Saari has been able to sidestep many of the prohibitively difficult combinatorial obstructions that are often associated with voting analysis.   In doing so, he draws our attention to the specialized roles played by a handful of specific subspaces of the vector space of profiles.  This geometric approach has led to an impressive number of unexpected results in voting theory, as well as a refined understanding of many fundamental and important results in the field (see \cite{saari-2001-chaotic} and \cite{saari-2001-decisions} for engaging overviews of the subject).  

Ever since we began thinking about the mathematics of voting, we have been intrigued by the prominent role played by symmetry arguments in Saari's work (see \cite{saari-1988-symmetry}, \cite{saari-1999-explaining}, and \cite{saari-2001-chaotic} for gentle introductions to such arguments).  In particular, it seemed to us that many of the symmetry-based ideas we were encountering could be explained easily if we only had the right algebraic framework.  In this paper, we describe just such a framework.  More specifically, we show how voting may be viewed naturally from an algebraic perspective by viewing profiles as elements of certain well-studied $\mathbb{Q}S_n$-modules.  

By using only a handful of simple combinatorial objects (e.g., tabloids) and some basic ideas from representation theory (e.g., Schur's Lemma), we are able to recast and extend some of Saari's well-known results.  For example, we recover a result concerning the important relationship between the Borda count and pairwise voting when voters return full rankings of the candidates (Theorem~\ref{theorem:  borda full}).  We then extend this result to a situation in which voters return partial rankings of the candidates (Theorem~\ref{theorem:  borda analogue}).  In the process, we construct an infinite  family of ``Borda-like" voting procedures.  With the help of our main theorem (Theorem~\ref{theorem:  big one}), we also, for example, address the relationship between positional voting and approval voting (Theorem~\ref{theorem:  approval voting}).

Our experience to date has convinced us that approaching the study of voting from an algebraic perspective can be incredibly illuminating when it comes to understanding the mathematical underpinnings of many different voting procedures.  In fact, we see this paper as a first step toward what might eventually be called \emph{algebraic voting theory}.   Although the ideas presented here are just the tip of the iceberg, we believe they will be of great interest to voting theorists and enthusiasts alike.


\section{Voting on Tabloids}

We begin by introducing combinatorial objects called \emph{tabloids}.  These objects play an important role in the representation theory of the symmetric group (see, for example, \cite{sagan-2001-symmetric}).  Tabloids also appear in the analysis of \emph{partially ranked data}, which includes the type of voting data we will be considering throughout this paper (see \cite{diaconis-1988-group} and \cite{marden-1995-analyzing}).  

Let $n$ be a positive integer.  A \emph{composition} of $n$ is a sequence $\lambda = (\lambda_1,\dots,\lambda_m)$ of positive integers whose sum is $n$.  If $\lambda_1 \ge \cdots \ge \lambda_m$, then $\lambda$ is a \emph{partition} of $n$.  For example, $\lambda = (2, 1, 1, 3)$ is a composition of $7$, but not a partition of $7$, and $(4, 2, 1, 1)$ is a partition (and therefore also a composition) of $8$.  

\begin{figure}[h]
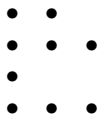

\begin{center}
\begin{tabular}{ccc} 
$\bullet$  & $\bullet$  \\
$\bullet$  & $\bullet$  & $\bullet$   \\
$\bullet$  \\
$\bullet$  & $\bullet$  & $\bullet$   \\
\end{tabular}
\end{center}
\caption{The Ferrers diagram of shape $(2,3,1,3)$.}
\label{Ferrers diagram}
\end{figure}

The \emph{Ferrers diagram of shape $\lambda$} is the left-justified array of dots with $\lambda_i$ dots in the $i$th row (see Figure \ref{Ferrers diagram}).  If the dots of a Ferrers diagram of shape $\lambda$ are replaced by boxes containing the numbers $1,\dots,n$ without repetition, then we create a \emph{Young tableau of shape $\lambda$}.  Two Young tableaux are said to be \emph{row equivalent} if they differ only by a permutation of the entries within the rows of each tableau.  An equivalence class of tableaux under this relation is called a \emph{tabloid of shape $\lambda$}.

We will denote a tabloid by first forming a representative tableau and then removing the vertical dividers within each row (see Figure \ref{tabloid}).  For convenience, we will usually choose the representative tableau whose entries in each row are in ascending order.  

\begin{figure}[h]
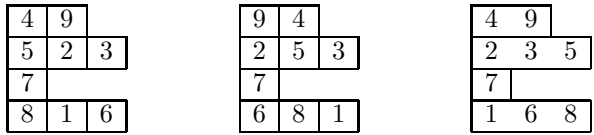
\label{tabloid}
\begin{center}
\begin{tabular}{| c | c | c |} \cline{1-2}
4 & 9 & \multicolumn{1}{| c}{ } \\ \cline{1-3}  
5 & 2 & 3 \\ \cline{1-3}
7 & \multicolumn{2}{| c}{} \\ \cline{1-3}
8 & 1 & 6 \\ \cline{1-3} 
\end{tabular}
\hspace{.5in}
\begin{tabular}{| c | c | c |} \cline{1-2}
9 & 4 & \multicolumn{1}{| c}{ } \\ \cline{1-3}  
2 & 5 & 3 \\ \cline{1-3}
7 & \multicolumn{2}{| c}{} \\ \cline{1-3}
6 & 8 & 1 \\ \cline{1-3} 
\end{tabular}
\hspace{.5in}
\begin{tabular}{| c c c |} \cline{1-2}
4 & 9 & \multicolumn{1}{| c}{ } \\ \cline{1-3}  
2 & 3 & 5 \\ \cline{1-3}
7 & \multicolumn{2}{| c}{} \\ \cline{1-3}
1 & 6 & 8 \\ \cline{1-3} 
\end{tabular}
\end{center}
\caption{Two equivalent tableaux and their tabloid.}
\label{tabloid}
\end{figure}

Let $X^\lambda$ denote the set of tabloids of shape $\lambda$.  Tabloids can be used to index voting data.  For example, suppose there are $n$ candidates, $ c_1, \dots, c_n$, in an election.  If the voters have been asked to vote by returning a list of the candidates in order of preference, from most to least favored, then each voter is essentially being asked to choose a tabloid from the set $X^{(1, \dots, 1)}$.  For example, if there are $n = 3$ candidates, then each voter must return one of the following tabloids from $X^{(1,1,1)}$:
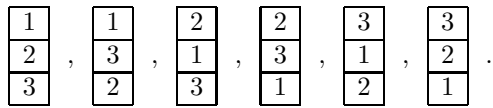
\begin{center}
\begin{tabular}{| c |} \hline
1 \\ \hline
2 \\ \hline
3 \\ \hline
\end{tabular} \hspace{.05in},\hspace{.05in}
\begin{tabular}{| c |} \hline
1 \\ \hline
3 \\ \hline
2 \\ \hline
\end{tabular} \hspace{.05in},\hspace{.05in}
\begin{tabular}{| c |} \hline
2 \\ \hline
1 \\ \hline
3 \\ \hline
\end{tabular} \hspace{.05in},\hspace{.05in}
\begin{tabular}{| c |} \hline
2 \\ \hline
3 \\ \hline
1 \\ \hline
\end{tabular} \hspace{.05in},\hspace{.05in}
\begin{tabular}{| c |} \hline
3 \\ \hline
1 \\ \hline
2 \\ \hline
\end{tabular} \hspace{.05in},\hspace{.05in}
\begin{tabular}{| c |} \hline
3 \\ \hline
2 \\ \hline
1 \\ \hline
\end{tabular} \ .
\end{center}
In this case, a voter who returns the second tabloid above prefers $c_1$ to $c_3$ to $c_2$, whereas a voter who returns the last tabloid prefers $c_3$ to $c_2$ to $c_1$. 

On the other hand, suppose the voters are asked to simply vote for their favorite candidate.  We could certainly obtain this information from the choices they made from $X^{(1, \dots, 1)}$ by focusing only on top-ranked candidates.  It may be much easier (for us and them), however, to have them choose from the set $X^{(1, n-1)}$.  For example, if $n = 4$, then we would be asking our voters to choose one of the following tabloids from $X^{(1,3)}$: 

\begin{center}
\begin{tabular}{| c c c |} \cline{1-1}
1 & \multicolumn{2}{| c}{} \\ \cline{1-3}
2 & 3 & 4  \\ \cline{1-3}
\end{tabular} \hspace{.05in},\hspace{.05in}
\begin{tabular}{| c c c |} \cline{1-1}
2 & \multicolumn{2}{| c}{} \\ \cline{1-3}
1 & 3 & 4  \\ \cline{1-3}
\end{tabular} \hspace{.05in},\hspace{.05in}
\begin{tabular}{| c c c |} \cline{1-1}
3 & \multicolumn{2}{| c}{} \\ \cline{1-3}
1 & 2 & 4  \\ \cline{1-3}
\end{tabular} \hspace{.05in},\hspace{.05in}
\begin{tabular}{| c c c |} \cline{1-1}
4 & \multicolumn{2}{| c}{} \\ \cline{1-3}
1 & 2 & 3  \\ \cline{1-3}
\end{tabular} \ .
\end{center}

In this case, by choosing the third tabloid above, a voter is saying that her favorite candidate is $c_3$ (and that she is indifferent, as far as her vote is concerned, to candidates $c_1$, $c_2$, and $c_4$).  

When voters are asked to provide their rankings of the candidates by choosing a tabloid from $X^{(1, \dots, 1)}$, we say that they are giving \emph{full rankings} of the candidates; if they are choosing tabloids from $X^\lambda$ where $\lambda \ne (1, \dots, 1)$, then we say that they are giving \emph{partial rankings} of the candidates.  

To determine the winner of an election based on choosing tabloids from $X^\lambda$, we need to know the number of votes received by each $x \in X^\lambda$.  With that in mind, let $\mathbf{p}:  X^\lambda \to \mathbb{N}$ be the function such that $\mathbf{p}(x)$ is the number of voters that voted for the tabloid $x$.  The function $\mathbf{p}$ is called a \emph{profile}.  

In fact, we will call every function $\mathbf{p}:  X^\lambda \to \mathbb{Q}$ a profile.  We do this primarily because, unlike functions from $X^\lambda$ to $\mathbb{N}$, the set $M^\lambda = \{ \mathbf{p}:  X^\lambda \to \mathbb{Q} \}$ forms a vector space (over $\mathbb{Q}$) where, for all $\mathbf{p}, \mathbf{q} \in M^\lambda$ and all $\alpha \in \mathbb{Q}$, $(\mathbf{p} + \mathbf{q})(x) = \mathbf{p}(x) + \mathbf{q}(x)$ and $(\alpha \mathbf{p})(x) = \alpha \mathbf{p}(x)$.  Note that the dimension of $M^\lambda$ is $| X^\lambda |$, and that the \emph{indicator functions} form a basis for $M^\lambda$ (where the indicator function for $x \in X^\lambda$ is the function that is $1$ on $x$ and $0$ on every other tabloid in $X^\lambda$).  

When it is convenient, we will also view a profile $\mathbf{p} \in M^\lambda$ as a formal linear combination of the tabloids in $X^\lambda$, where the coefficient in front of the tabloid $x$ is $\mathbf{p}(x)$.  For example, if $n = 3$ and there are eleven voters, then our profile (which is an example used in Chp.~2 of \cite{saari-2001-chaotic}) might look like

\begin{center}
3 \
\begin{tabular}{| c |} \hline
1 \\ \hline
2 \\ \hline
3 \\ \hline
\end{tabular}
+ 2 \
\begin{tabular}{| c |} \hline
1 \\ \hline
3 \\ \hline
2 \\ \hline
\end{tabular}
+ 0 \
\begin{tabular}{| c |} \hline
2 \\ \hline
1 \\ \hline
3 \\ \hline
\end{tabular}
+ 2 \
\begin{tabular}{| c |} \hline
2 \\ \hline
3 \\ \hline
1 \\ \hline
\end{tabular} 
+ 0 \
\begin{tabular}{| c |} \hline
3 \\ \hline
1 \\ \hline
2 \\ \hline
\end{tabular} 
+ 4 \
\begin{tabular}{| c |} \hline
3 \\ \hline
2 \\ \hline
1 \\ \hline
\end{tabular}
\end{center}
where three voters chose the first list, two voters chose the second list, zero voters chose the third list, and so on.  

Although viewing profiles as formal linear combinations of tabloids can be useful at times, we will typically view profiles as column vectors in $\mathbb{Q}^{|X^\lambda|}$.  We do this by choosing the indicator functions of $M^\lambda$ as a basis, which we order based on the lexicographic ordering of the tabloids.  For example, we would encode the above profile as the vector
\[
\mathbf{p} = 
\begin{bmatrix}
3 \\
2 \\
0 \\
2 \\ 
0 \\
4 \\
\end{bmatrix}
\begin{matrix}
123 \\
132 \\
213 \\
231 \\
312 \\
321 \\
\end{matrix}
\]
where we have placed labels corresponding to the rankings to the right of the vector so that the correspondence can be seen more easily.


\section{Positional Voting Procedures}

Now that we know what profiles look like, we can turn our attention to voting procedures.  In this paper, we will primarily focus our attention on voting procedures that assign points to each candidate based on their position in a voter's choice of a tabloid.  

Let $\lambda = (\lambda_1, \dots, \lambda_m)$ be a composition of $n$, let $\mathbf{w} = [w_1, \dots, w_m]^t$ be a vector in $\mathbb{Q}^m$, and suppose we are given a profile $\mathbf{p} \in M^\lambda$.  We use the vector $\mathbf{w}$, which is called a \emph{weighting vector}, to assign points to each candidate, and the candidate who receives the most points is declared to be the winner.  We do this as follows.  For each tabloid $x$, if candidate $c_i$ is in row $j$ of $x$, then she will be given $w_j \mathbf{p}(x)$ points.  Summing over all tabloids $x \in X^\lambda$ then determines the total number of points assigned to candidate $c_i$.  We will refer to this as a \emph{positional voting} procedure of type $\lambda$ based on $\mathbf{w}$.  

For example, consider our three-candidate and eleven-voter example above.  If we use a positional voting procedure with weighting vector $\mathbf{w} = [1, s, 0]^t$, where $0 \le s \le 1$, then candidate $c_1$ will receive $(3 \times 1) + (0 \times s) + (6 \times 0) = 5$ points, candidate $c_2$ will receive $(2 \times 1) + (7 \times s) + (2 \times 0) = 2 + 7s$ points, and candidate $c_3$ will receive $(4 \times 1) + (4 \times s) + (3 \times 0) = 4 + 4s$ points.  We encode these points in the \emph{results vector} $\mathbf{r} = [5, 2 + 7s, 4 + 4s]^t$.  

The point, of course, of using the parameter $s$ in this example is that it highlights one of the things that makes voting theory so interesting:  for a fixed profile $\mathbf{p}$, the outcomes of an election can vary wildly with the choice of an election procedure.  In fact, as Saari points out in \cite{saari-2001-chaotic},
\begin{quote}
``...rather than reflecting the views of the voters, it is entirely possible for an election outcome to more accurately reflect the choice of an election procedure."
\end{quote}
To see this, note that in the above example, when $s = 0$, we have the well-known plurality voting procedure (``vote for your favorite"), in which case $c_1$ wins with a results vector of $\mathbf{r} = [5, 2, 4]^t$.  When $s = 1$ (``vote for everyone but your least favorite"), $c_2$ wins with a results vector of $\mathbf{r} = [5, 9, 8]^t$, and when $s = 1/2$, $c_3$ wins with a results vector of $\mathbf{r} = [5, 5 \frac{1}{2}, 6]^t$.  (See Chp.~2 of \cite{saari-2001-chaotic} for this and other \emph{voting paradoxes}.)

One attractive feature of positional voting is that the results vectors they produce can be realized as the product of a matrix, which we will denote by $T_{\mathbf{w}}$, and the profile in question.  For example, if $\mathbf{p}$ and $\mathbf{w} = [1, s, 0]^t$ are as above, then 
\begin{equation*}
T_{\mathbf{w}} (\mathbf{p}) 
= 
\begin{bmatrix}
1 & 1 & s & 0 & s & 0 \\
s & 0 & 1 & 1 & 0 & s \\
0 & s & 0 & s & 1 & 1
\end{bmatrix}
\begin{bmatrix}
3 \\ 2 \\ 0 \\ 2 \\ 0 \\ 4
\end{bmatrix}
=
\begin{bmatrix}
5 \\ 2 + 7s \\ 4 + 4s
\end{bmatrix}
=
\mathbf{r}.
\end{equation*}
In fact, and this is a simple but key insight, every positional voting method of type $\lambda$ based on the weighting vector $\mathbf{w}$ can be viewed as a linear transformation $T_{\mathbf{w}}:  M^\lambda \to M^{(1, n-1)}$ since we may use the tabloids in $X^{(1, n-1)}$ to index the set of candidates.  For example, we have $T_{[1, s, 0]^t}:  M^{(1,1,1)} \to M^{(1, 2)}$ in the example above. 

Recognizing that $T_{\mathbf{w}}:  M^\lambda \to M^{(1, n-1)}$ is a linear transformation is certainly useful, but there are other, hidden algebraic structures within the framework we have just constructed.  To begin to see some of these structures, let $\lambda = (\lambda_1, \dots, \lambda_m)$ be a composition of $n$.  There are $n!$ tableaux of shape $\lambda$.  In fact, we may think of a tableau as a full ranking of the candidates by reading, left to right, top to bottom, the entries of the tableau.  For example, the tableau
\begin{center}
\begin{tabular}{| c | c |} \cline{1-2}
2 & 5 \\ \cline{1-2}
1 & 3 \\ \cline{1-2}
4 & \multicolumn{1}{| c}{} \\ \cline{1-1} 
\end{tabular}
\end{center}
would correspond to the full ranking
\begin{center}
\begin{tabular}{| c |} \cline{1-1}
2 \\ \cline{1-1}
5 \\ \cline{1-1}
1 \\ \cline{1-1}
3 \\ \cline{1-1}
4 \\ \cline{1-1}
\end{tabular} \ .
\end{center}
This is useful, because we may then view a vote for the tabloid
\begin{center}
\begin{tabular}{| c c |} \cline{1-2}
2 & 5 \\ \cline{1-2}
1 & 3 \\ \cline{1-2}
4 & \multicolumn{1}{| c}{} \\ \cline{1-1} 
\end{tabular}
\end{center}
as coming from anyone whose full ranking of the candidates corresponds to one of the tableaux in the equivalence class of the tabloid.  For example, the above tabloid could correspond to any of the following full rankings:
\begin{center}
\begin{tabular}{| c |} \cline{1-1}
2 \\ \cline{1-1}
5 \\ \cline{1-1}
1 \\ \cline{1-1}
3 \\ \cline{1-1}
4 \\ \cline{1-1}
\end{tabular} \hspace{.05in},\hspace{.05in}
\begin{tabular}{| c |} \cline{1-1}
2 \\ \cline{1-1}
5 \\ \cline{1-1}
3 \\ \cline{1-1}
1 \\ \cline{1-1}
4 \\ \cline{1-1}
\end{tabular} \hspace{.05in},\hspace{.05in}
\begin{tabular}{| c |} \cline{1-1}
5 \\ \cline{1-1}
2 \\ \cline{1-1}
1 \\ \cline{1-1}
3 \\ \cline{1-1}
4 \\ \cline{1-1}
\end{tabular} \hspace{.05in},\hspace{.05in}
\begin{tabular}{| c |} \cline{1-1}
5 \\ \cline{1-1}
2 \\ \cline{1-1}
3 \\ \cline{1-1}
1 \\ \cline{1-1}
4 \\ \cline{1-1}
\end{tabular} \ .
\end{center}

This simple insight allows us to view profiles $\mathbf{p} \in M^\lambda$ as elements of $M^{(1, \dots, 1)}$ in the following manner.  Let $f$ be the number of full rankings corresponding to a tabloid of shape $\lambda$.  Define $\iota:  M^\lambda \to M^{(1, \dots, 1)}$ (think ``inclusion") by mapping each tabloid (i.e., each indicator function) to the sum of its corresponding full rankings times $1/f$.  For example, the tabloid of shape $\lambda = (2, 2, 1)$ above would be mapped to 
\[
\frac{1}{4}
\left(
\text{
\begin{tabular}{| c |} \cline{1-1}
2 \\ \cline{1-1}
5 \\ \cline{1-1}
1 \\ \cline{1-1}
3 \\ \cline{1-1}
4 \\ \cline{1-1}
\end{tabular}
}
+
\text{
\begin{tabular}{| c |} \cline{1-1}
2 \\ \cline{1-1}
5 \\ \cline{1-1}
3 \\ \cline{1-1}
1 \\ \cline{1-1}
4 \\ \cline{1-1}
\end{tabular}
}
+
\text{
\begin{tabular}{| c |} \cline{1-1}
5 \\ \cline{1-1}
2 \\ \cline{1-1}
1 \\ \cline{1-1}
3 \\ \cline{1-1}
4 \\ \cline{1-1}
\end{tabular} 
}
+
\text{
\begin{tabular}{| c |} \cline{1-1}
5 \\ \cline{1-1}
2 \\ \cline{1-1}
3 \\ \cline{1-1}
1 \\ \cline{1-1}
4 \\ \cline{1-1}
\end{tabular}
}
\right).
\]
If we define $\pi:  M^{(1, \dots, 1)} \to M^\lambda$ (think ``projection") by mapping each full ranking to the tabloid of shape $\lambda$ that contains it, then note that $\pi \circ \iota:  M^\lambda \to M^\lambda$ is the identity transformation.

Let $\mathbf{p} \in M^\lambda$, and define $\overline{\mathbf{p}} = \iota(\mathbf{p})$.  In other words, each profile $\mathbf{p} \in M^\lambda$ maybe viewed as a profile $\overline{\mathbf{p}} \in M^{(1, \dots, 1)}$ that is constant on the equivalence classes that form the tabloids of shape $\lambda$.  Moreover, if $\mathbf{w} = [w_1, \dots, w_m]^t$ is a weighting vector associated to $\lambda$, and we define $\overline{\mathbf{w}}$ to be the weighting vector in $\mathbb{Q}^n$ whose first $\lambda_1$ entries are equal to $w_1$, whose next $\lambda_2$ entries are equal to $w_2$, and so on, then
\[
T_{\mathbf{w}}(\mathbf{p}) = T_{\overline{\mathbf{w}}}(\overline{\mathbf{p}}).
\]
For example, suppose $\lambda = (2, 1)$, $\mathbf{w} = [3, 0]^t$, and 
\begin{center}
$\mathbf{p} = $
5 \
\begin{tabular}{| c c |} \cline{1-2}
1 & 2 \\ \cline{1-2}
3 & \multicolumn{1}{| c}{} \\ \cline{1-1} 
\end{tabular}
+ 4 \
\begin{tabular}{| c c |} \cline{1-2}
1 & 3 \\ \cline{1-2}
2 & \multicolumn{1}{| c}{} \\ \cline{1-1} 
\end{tabular}
+ 7 \
\begin{tabular}{| c c |} \cline{1-2}
2 & 3 \\ \cline{1-2}
1 & \multicolumn{1}{| c}{} \\ \cline{1-1} 
\end{tabular}.
\end{center}
Then $\overline{\mathbf{w}} = [3,3,0]^t$, 
\[
\overline{\mathbf{p}} = 
\frac{5}{2} 
\left(
\text{ 
\begin{tabular}{| c |} \cline{1-1}
1 \\ \cline{1-1}
2 \\ \cline{1-1}
3 \\ \cline{1-1} 
\end{tabular}
\ + 
\begin{tabular}{| c |} \cline{1-1}
2 \\ \cline{1-1}
1 \\ \cline{1-1}
3 \\ \cline{1-1} 
\end{tabular}
}
\right)
+
2
\left(
\text{ 
\begin{tabular}{| c |} \cline{1-1}
1 \\ \cline{1-1}
3 \\ \cline{1-1}
2 \\ \cline{1-1} 
\end{tabular}
\ + 
\begin{tabular}{| c |} \cline{1-1}
3 \\ \cline{1-1}
1 \\ \cline{1-1}
2 \\ \cline{1-1} 
\end{tabular}
}
\right)
+
\frac{7}{2} 
\left(
\text{ 
\begin{tabular}{| c |} \cline{1-1}
2 \\ \cline{1-1}
3 \\ \cline{1-1}
1 \\ \cline{1-1} 
\end{tabular}
\ + 
\begin{tabular}{| c |} \cline{1-1}
3 \\ \cline{1-1}
2 \\ \cline{1-1}
1 \\ \cline{1-1} 
\end{tabular}
}
\right),
\]
and
\begin{equation*}
T_{\mathbf{w}} (\mathbf{p}) 
= 
\begin{bmatrix}
3 & 3 & 0 \\
3 & 0 & 3 \\
0 & 3 & 3 \\
\end{bmatrix}
\begin{bmatrix}
5 \\ 4 \\ 7
\end{bmatrix}
=
\begin{bmatrix}
3 & 3 & 3 & 0 & 3 & 0 \\
3 & 0 & 3 & 3 & 0 & 3 \\
0 & 3 & 0 & 3 & 3 & 3
\end{bmatrix}
\begin{bmatrix}
\frac{5}{2} \\ 2 \\ \frac{5}{2} \\ \frac{7}{2} \\ 2 \\ \frac{7}{2}
\end{bmatrix}
=
T_{\overline{\mathbf{w}}}(\overline{\mathbf{p}}).
\end{equation*}

\section{Representation Theory and the Symmetric Group}

We have just seen that profiles can be viewed as vectors, and that each positional voting procedure can be viewed as a linear transformation.  In this section, we describe the role that the symmetric group plays in positional voting.   In particular, as we illustrate below, each map $T_\mathbf{w}:  M^\lambda \to M^{(1, n-1)}$ is more than just a linear transformation---it is a \emph{$\mathbb{Q}S_n$-module homomorphism}.  

To explain, let $\lambda$ be a composition of $n$.  The symmetric group $S_n$ acts naturally on the set $X^\lambda$ of tabloids of shape $\lambda$ by permuting the entries of the tabloids.  For example, if $n = 5$ and $\sigma = (1 \ 3) (2 \ 5 \ 4)$, then 
\begin{center}
$\sigma$
\begin{tabular}{| c c c |} \cline{1-3}
2 & 3 & 5 \\ \cline{1-3}
1 & 4 & \multicolumn{1}{| c}{} \\ \cline{1-2}
\end{tabular} 
\hspace{.05in} $=$ \hspace{.05in}
\begin{tabular}{| c c c |} \cline{1-3}
$\sigma(2)$ & $\sigma(3)$ & $\sigma(5)$ \\ \cline{1-3}
$\sigma(1)$ & $\sigma(4)$ & \multicolumn{1}{| c}{} \\ \cline{1-2}
\end{tabular}
\hspace{.05in} $=$ \hspace{.05in}
\begin{tabular}{| c c c |} \cline{1-3}
5 & 1 & 4 \\ \cline{1-3}
3 & 2 & \multicolumn{1}{| c}{} \\ \cline{1-2}
\end{tabular}.
\end{center}
We may extend the action of $S_n$ on $X^\lambda$ to an action of $S_n$ on the profile space $M^\lambda$ by defining $(\sigma \mathbf{p})(x) = \mathbf{p}(\sigma^{-1}x)$.  Moreover, we may extend this action to an action of the group ring $\mathbb{Q}S_n$ on $M^\lambda$ where if $a = \sum_{\sigma \in S_n} a_\sigma \sigma$, then $(a \mathbf{p})(x) = \sum_{\sigma \in S_n} a_\sigma \mathbf{p}(\sigma^{-1}x)$.  

The action of $\mathbb{Q}S_n$ on $M^\lambda$ turns $M^\lambda$ into something called a \emph{$\mathbb{Q}S_n$-module}.  We will not discuss modules in general in this paper (see, for example, \cite{dummit-foote-2004-abstract} for a nice introduction to modules).  We will, however, describe some of the more useful implications of this realization as far as voting is concerned.  

To begin, note that, for each element $a \in \mathbb{Q}S_n$, there is a linear transformation $L_a: M^\lambda \to M^\lambda$ defined by setting $L_a(\mathbf{p}) = a\mathbf{p}$.  If we let $\text{End}(M^\lambda)$ denote the ring of linear transformations from $M^\lambda$ to itself (the ring of \emph{endomorphisms}), then this defines a ring homomorphism $\rho:  \mathbb{Q}S_n \to \text{End}(M^\lambda)$.  

The homomorphism $\rho$ is an example of a \emph{representation} of $\mathbb{Q}S_n$ since each element $a \in \mathbb{Q}S_n$ may be ``represented" by a linear transformation $\rho(a) = L_a \in \text{End}(M^\lambda)$.  Furthermore, if we restrict $\rho$ to the elements of $S_n$, then the images are all invertible linear transformations.  In other words, we have a group homomorphism $\rho \downarrow_{S_n} :  S_n \to \text{Aut}(M^\lambda) = \text{GL}(M^\lambda)$, which is a \emph{representation} of $S_n$.  

The fact that $M^\lambda$ is a $\mathbb{Q}S_n$-module tells us that we can write $M^\lambda$ as a direct sum of \emph{$\mathbb{Q}S_n$-submodules}, which are subspaces that are invariant under the action of $\mathbb{Q}S_n$ (and are therefore $\mathbb{Q}S_n$-modules themselves).  In other words, there are decompositions
\[
M^\lambda = M_1 \oplus \cdots \oplus M_k
\]
such that, for all $a \in \mathbb{Q}S_n$ and $\mathbf{p}_i \in M_i$, $a \mathbf{p}_i \in M_i$.  

A submodule $U$ of a module $M$ is said to be \emph{simple} if $U \ne 0$, and the only submodules of $U$ are $0$ and $U$ itself.  It turns out that, up to isomorphism, there are only a finite number of distinct simple $\mathbb{Q}S_n$-modules.  Furthermore, these simple modules can be (and typically are) indexed by the partitions of $n$.  In this paper, we use a well-known indexing scheme (see \cite{sagan-2001-symmetric}), and we denote the simple $\mathbb{Q}S_n$-module indexed by the partition $\mu$ of $n$ by $S^\mu$.  

As an example, it turns out that $M^{(1, 2)} \cong S^{(3)} \oplus S^{(2,1)}$ and $M^{(1,1,1)} \cong S^{(3)} \oplus S^{(2,1)} \oplus S^{(2,1)} \oplus S^{(1,1,1)}$.  To make these decompositions more concrete, consider the following decomposition of $M^{(1,2)}$:
\begin{equation*}
M^{(1,2)} = 
\langle 
\begin{bmatrix} 
1 \\ 1 \\ 1
\end{bmatrix}
\rangle
\oplus
\langle 
\begin{bmatrix} 
1 \\ -1 \\ 0
\end{bmatrix},
\begin{bmatrix} 
1 \\ 0 \\ -1
\end{bmatrix}
\rangle.
\end{equation*}

Both of the subspaces on the right are invariant under the action of $\mathbb{Q}S_3$.  The first space, the span of the ``all-ones" vector, is obviously invariant under the action of $\mathbb{Q}S_3$.  After all, the entries in these vectors are all equal, and permuting the entries of such vectors does not change this fact.  

The second subspace, which is the orthogonal complement (with respect to the usual dot product) of the first subspace, is the subspace of vectors whose entries sum to zero.  Since this property too is preserved under the action of $\mathbb{Q}S_3$, it too is a submodule of $M^{(1,2)}$.  

Both of these subspaces are simple $\mathbb{Q}S_3$-modules.  The first space is isomorphic to $S^{(3)}$, and in general, we use $S^{(n)}$ to denote this ``trivial", one-dimensional $\mathbb{Q}S_n$-module that has the property that every element of $S_n$ acts as the identity linear transformation.  The second space is isomorphic to $S^{(2,1)}$, and in general, $S^{(n-1,1)}$ is the $(n-1)$-dimensional simple module that is isomorphic to the orthogonal complement of the span of the all-ones vector in the $n$-dimensional $\mathbb{Q}S_n$-module $M^{(1, n-1)} \cong M^{(n-1,1)}$.    

To see how all of this is related to voting, let $\lambda = (\lambda_1, \dots, \lambda_m)$ be a composition of $n$, let $\mathbf{w} = [w_1, \dots, w_m]^t$ be a weighting vector, and consider the positional voting procedure $T_\mathbf{w}:  M^\lambda \to M^{(1,n-1)}$ of shape $\lambda$ based on $\mathbf{w}$.  First, note that the linear transformation $T_\mathbf{w}$ is actually a \emph{$\mathbb{Q}S_n$-module homomorphism}.  In other words, if $a \in \mathbb{Q}S_n$, then $T_\mathbf{w}(a \mathbf{p}) = a T_\mathbf{w}(\mathbf{p})$.  After all, if $\sigma \in S_n$, then we are simply asking that $T_{\mathbf{w}}(\sigma \mathbf{p}) = \sigma T_{\mathbf{w}}(\mathbf{p})$, which is the same as saying that, ``if we permute the labels on the candidates, then we need only permute the original scores in the same way."  We are just extending this notion, which is called \emph{neutrality} (see, for example, \cite{saari-valognes-1998-geometry}), to the action of the entire group ring. 

Viewing $T_\mathbf{w}$ as a $\mathbb{Q}S_n$-module homomorphism is helpful for a couple of reasons.  First, if $T:  M \to N$ is a module homomorphism, then the kernel of $T$ is a submodule of $M$, and the image of $T$ is a submodule of $N$.  In fact, if we define the \emph{effective space} $E(T)$ of $T$ to be the orthogonal complement (with respect to the usual dot product) of the kernel of $T$, then $E(T) \cong T(M)$ as modules.  Therefore, knowing something about $E(T_\mathbf{w})$ might help us say something about $T_\mathbf{w}$.  

Second, once we know we are dealing with module homomorphisms, we can look to an elementary but immensely useful theorem for insight: 
\begin{quote}
\textbf{Schur's Lemma}.  Any nonzero module homomorphism between simple modules is an isomorphism.
\end{quote}  
How can we make use of Schur's Lemma in our study of voting?  Consider our positional voting procedure $T_\mathbf{w} : M^\lambda \to M^{(1, n-1)}$.  As we noted above, the module $M^{(1,n-1)}$ is isomorphic to a direct sum of the simple modules $S^{(n)}$ and $S^{(n-1,1)}$.  This means that any simple submodule $U$ of $M^\lambda$ that is not isomorphic to $S^{(n)}$ or $S^{(n-1,1)}$ \emph{must} be in the kernel of $T_\mathbf{w}$, i.e., such a submodule only contains information that will have absolutely no effect on the results of the election. 

For example, if $n = 3$ and we have $T_\mathbf{w} :  M^{(1,1,1)} \to M^{(1,2)}$, then since $M^{(1, 2)} \cong S^{(3)} \oplus S^{(2,1)}$ and $M^{(1,1,1)} \cong S^{(3)} \oplus S^{(2,1)} \oplus S^{(2,1)} \oplus S^{(1,1,1)}$, we know that the kernel  of $T_\mathbf{w}$ must contain exactly one copy of $S^{(1,1,1)}$ and at least one copy of $S^{(2,1)}$.  We will say more about this below.  Before we do, however, we need one more idea.  Although it may appear simple on the surface, it will play a major role in what follows. 

First, note that there is a bijection between the full rankings in $X^{(1,\dots,1)}$ and the permutations in $S_n$, where the tabloid
\begin{center}
\begin{tabular}{| c |} \cline{1-1}
$i_1$ \\ \cline{1-1}
$i_2$ \\ \cline{1-1}
$\vdots$ \\ \cline{1-1} 
$i_{n-1}$ \\ \cline{1-1}
$i_{n}$ \\ \cline{1-1}
\end{tabular}
\end{center}
is mapped to the permutation $\sigma$ with the property that $\sigma(j) = i_j$.  It follows that we may view each profile $\mathbf{p} \in M^{(1, \dots, 1)}$ as an element of $\mathbb{Q}S_n$, where we simply replace each full ranking with its associated permutation.  This means that if $\mathbf{p} \in M^\lambda$, then we may view $\overline{\mathbf{p}} \in M^{(1, \dots, 1)}$ as an element of $\mathbb{Q}S_n$ that is constant on the left cosets of the subgroup of $S_n$ that fixes the tabloid containing the tableau corresponding to the identity of $S_n$.  (We will use this fact in the proof of Theorem~\ref{theorem:  big one}.)

How is viewing a profile as an element of $\mathbb{Q}S_n$ helpful?  It means that our positional voting procedures are more than just linear transformations, and more than just $\mathbb{Q}S_n$-module homomorphisms---they are the results of \emph{profiles acting on weighting vectors}!  More specifically, if $\mathbf{p} \in \mathbb{Q}S_n$ is a profile, and $\mathbf{w} \in \mathbb{Q}^n \cong M^{(1, n-1)}$ is a weighting vector, then 
\[
T_{\mathbf{w}}(\mathbf{p}) = \mathbf{p} \mathbf{w}.
\]

For example, in the three-candidate and eleven-voter example introduced above, note that if $e \in S_3$ is the identity and we use the usual cycle notation for the other elements in $S_3$, then 
\begin{align*}
T_{\mathbf{w}} (\mathbf{p}) &= 
\begin{bmatrix}
1 & 1 & s & 0 & s & 0 \\
s & 0 & 1 & 1 & 0 & s \\
0 & s & 0 & s & 1 & 1
\end{bmatrix}
\begin{bmatrix}
3 \\ 2 \\ 0 \\ 2 \\ 0 \\ 4
\end{bmatrix} \\
&= 
3
\begin{bmatrix}
1 \\ s \\ 0
\end{bmatrix}
+
2
\begin{bmatrix}
1 \\ 0 \\ s
\end{bmatrix}
+
2
\begin{bmatrix}
0 \\ 1 \\s 
\end{bmatrix}
+
4
\begin{bmatrix}
0 \\ s \\ 1
\end{bmatrix} \\
&=
\left( 3e + 2 (23) + 2 (123) + 4 (13) \right)
\begin{bmatrix}
1 \\ s \\ 0
\end{bmatrix} \\
&=
\mathbf{p} \mathbf{w}.
\end{align*}

This realization---that if $\mathbf{p} \in M^{(1, \dots, 1)}$, then $T_\mathbf{w}(\mathbf{p}) = \mathbf{p} \mathbf{w}$, i.e., ``the acted upon has become the actor"---leads us to our first theorem.  Before we state the theorem, though, we need a few observations about weighting vectors whose entries sum to zero.  

It is important to keep in mind the relatively simple structure of the space $M^{(1, n-1)}$ in which our weighting vector $\mathbf{w}$ resides (together, of course, with all of our results vectors).  In particular, since $M^{(1,n-1)} \cong S^{(n)} \oplus S^{(n-1,1)}$, we may write $\mathbf{w} = \mathbf{1}_\mathbf{w} + \widehat{\mathbf{w}}$, where $\mathbf{1}_\mathbf{w} \in S^{(n)}$ is the projection of $\mathbf{w}$ onto the all-ones vector $\mathbf{1}$, and $\widehat{\mathbf{w}} \in S^{(n-1,1)}$ is the projection of $\mathbf{w}$ into the orthogonal $(n-1)$-dimensional subspace of vectors whose entries sum to zero.  For example, if $n = 3$ and $\mathbf{w} = [1, s, 0]^t$, then $\mathbf{1}_\mathbf{w} = [\frac{1+s}{3}, \frac{1+s}{3}, \frac{1+s}{3}]^t$ and $\widehat{\mathbf{w}} = [\frac{2-s}{3}, \frac{2s-1}{3}, \frac{-1-s}{3}]^t$.

Since
\[
T_\mathbf{w}(\mathbf{p}) = \mathbf{p} \mathbf{w} = \mathbf{p} (\mathbf{1}_\mathbf{w} + \widehat{\mathbf{w}})
= \mathbf{p} \mathbf{1}_\mathbf{w} + \mathbf{p} \widehat{\mathbf{w}}
\]
and $\mathbf{p} \mathbf{1}_\mathbf{w} \in S^{(n)}$, all of the information that will determine the outcome of the election is contained in the summand $\mathbf{p} \widehat{\mathbf{w}}$.  After all, $\mathbf{p} \mathbf{1}_\mathbf{w}$ is simply some multiple of the all-ones vector and will therefore not differentiate between any of the candidates.  Because of this, we will focus most of our attention on weighting vectors $\mathbf{w} \in M^{(1,n-1)} = \mathbb{Q}^n$ such that $\mathbf{w} = \widehat{\mathbf{w}}$, i.e., weighting vectors whose entries sum to zero.  

For convenience, we say that a vector in $\mathbb{Q}^n$ whose entries sum to zero is a \emph{sum-zero vector}.  We also say that such a vector is \emph{nontrivial} if it does not equal the zero vector (since the weighting vector $\mathbf{w} = \mathbf{0}$ would obviously lead to the trivial result that the candidates all tie with zero points).  A consequence of using a sum-zero weighting vector $\mathbf{w}$ is that the results vector $\mathbf{r} = T_\mathbf{w} (\mathbf{p})$ will also be a sum-zero vector.  Keep in mind, however, that the winner is still the candidate that receives the most number of points.  

\begin{theorem}\label{theorem:  big one}
Let $n \ge 2$, and let $\lambda = (\lambda_1, \dots, \lambda_m)$ be a partition of $n$.  Suppose that $\mathbf{w}_1, \dots, \mathbf{w}_k$ form a linearly independent set of weighting vectors in $\mathbb{Q}^m$ such that $\overline{\mathbf{w}_1}, \dots, \overline{\mathbf{w}_k}$ are sum-zero vectors.  If $\mathbf{r}_1, \dots, \mathbf{r}_k$ are any sum-zero results vectors in $\mathbb{Q}^n$, then there exist infinitely many profiles $\mathbf{p} \in M^\lambda$ such that $T_{\mathbf{w}_i}(\mathbf{p}) = \mathbf{r}_i$ for all $1 \le i \le k$.  
\end{theorem}

\begin{proof}
We first consider the full ranking case where $\lambda = (1, \dots, 1)$.  In this case, $\mathbf{w}_i = \overline{\mathbf{w}_i}$.  The sum-zero vectors in $\mathbb{Q}^n \cong M^{(1, n-1)}$ form a simple $\mathbb{Q}S_n$-submodule (that is isomorphic to $S^{(n-1,1)}$), which we will denote by $U$.  Since the weighting vectors are linearly independent, there exists a linear transformation $T: U \to U$ such that $T(\mathbf{w}_i) = \mathbf{r}_i$ for all $1 \le i \le k$.  

By a theorem of Burnside (see, for example, \cite{clausen-baum-1993-fast}), every linear transformation from a simple $\mathbb{Q}S_n$-module to itself can be realized as the action of some element in $\mathbb{Q}S_n$.  In other words, there is some $a \in \mathbb{Q}S_n$ such that $T(\mathbf{u}) = a \mathbf{u}$ for all $\mathbf{u} \in U$.  Moreover, there exist infinitely many $b \in \mathbb{Q}S_n$ such that $b \mathbf{u} = \mathbf{0}$ for all $\mathbf{u} \in U$.  If we set $\mathbf{p} = a + b$, then the theorem follows since $T_{\mathbf{w}_i}(\mathbf{p}) = \mathbf{p} \mathbf{w}_i = a \mathbf{w}_i + b \mathbf{w}_i = \mathbf{r}_i + \mathbf{0}$.  

For a general $\lambda$, note that if $\{ \mathbf{w}_1, \dots, \mathbf{w}_k \} \subset \mathbb{Q}^m$ is linearly independent, then so is $\{ \overline{\mathbf{w}_1}, \dots, \overline{\mathbf{w}_k} \} \subset \mathbb{Q}^n$.  We proved above that there exist infinitely many $\mathbf{p} \in \mathbb{Q}S_n$ such that $\mathbf{p} \overline{\mathbf{w}_i} = \mathbf{r}_i$.  But, as we saw earlier, we may replace $\mathbf{p}$ with $\mathbf{p}'$, where $\mathbf{p}'$ is constant on the equivalence classes that form the tabloids of shape $\lambda$.  We may therefore also view $\mathbf{p}'$ as an element of $M^\lambda$, in which case $T_{\mathbf{w}_i} (\mathbf{p}') = \mathbf{p} \overline{\mathbf{w}_i} = \mathbf{r}_i$ for all $1 \le i \le k$.  
\end{proof}

Theorem~\ref{theorem:  big one} is an extension of Saari's Theorem 1 in~\cite{saari-1984-ultimate} in two ways.  First, Saari's theorem is a statement about \emph{ordinal rankings}, i.e., the order in which the candidates finish in the election.  Our theorem says something about \emph{cardinal rankings},  i.e., the actual number of points that each candidate receives.  Second, whereas Saari's theorem focuses on the fully ranked situation, we address both the fully and partially ranked situations simultaneously.  

Both theorems essentially imply that as long as the weighting vectors $\mathbf{w}_1, \dots, \mathbf{w}_k$ are different enough, it could very well be the case that there is no relationship whatsoever between $T_{\mathbf{w}_1}(\mathbf{p}), \dots, T_{\mathbf{w}_k}(\mathbf{p})$.  Moreover, as Saari describes in \cite{saari-2001-chaotic} and \cite{saari-2001-decisions}, this is just the tip of the iceberg when it comes to answering the question, ``How bad can it get?"  Note, however, that our proof is decidedly algebraic in nature.  Whereas Saari's proof of Theorem 1 in~\cite{saari-1984-ultimate} uses facts about open mappings, our proof of Theorem~\ref{theorem:  big one} uses an important result by Burnside concerning the endomorphism ring of simple modules.


\section{Approval Voting}

Theorem~\ref{theorem:  big one} may be used to address paradoxical situations that arise in voting procedures related to positional voting.  For example, the \emph{approval voting} procedure asks a voter to return an (unordered) list of the candidates of whom she approves.  A candidate receives a point for each time she appears on such a list, and the candidate receiving the most points is declared the winner.  (For more on approval voting, see \cite{brams-fishburn-2007-approval}.)

We naturally assume that, if a voter were to return a fully ranked list of the candidates, then the candidates she would approve of would make up the top portion of her list.  We may therefore imagine a situation in which each voter is asked to return a fully ranked list of the candidates together with a cutoff point.  Candidates above the cutoff are those whom our voter approves of, and those below the cutoff are not.  

We will denote the cutoff point in a tableau with a blank space separating the ``approved" candidates from the other candidates.  For example, if $n = 3$, then the top row of Figure~\ref{figure:  tableaux with cutoffs} contains those tableaux with cutoffs that would be used by voters who only approve of their top candidate, whereas the second row contains those tableaux with cutoffs that would be used by those voters who approve of their top two candidates.  

\begin{figure}[h]
\begin{center}
\begin{tabular}{| c |} \hline
1 \\ \hline \hline
2 \\ \hline
3 \\ \hline
\end{tabular} \hspace{.2in}
\begin{tabular}{| c |} \hline
1 \\ \hline \hline
3 \\ \hline
2 \\ \hline
\end{tabular} \hspace{.2in}
\begin{tabular}{| c |} \hline
2 \\ \hline \hline
1 \\ \hline
3 \\ \hline
\end{tabular} \hspace{.2in}
\begin{tabular}{| c |} \hline
2 \\ \hline \hline
3 \\ \hline
1 \\ \hline
\end{tabular} \hspace{.2in}
\begin{tabular}{| c |} \hline
3 \\ \hline \hline
1 \\ \hline
2 \\ \hline
\end{tabular} \hspace{.2in}
\begin{tabular}{| c |} \hline
3 \\ \hline \hline
2 \\ \hline
1 \\ \hline
\end{tabular}
\end{center}

\bigskip

\begin{center}
\begin{tabular}{| c |} \hline
1 \\ \hline 
2 \\ \hline \hline 
3 \\ \hline
\end{tabular} \hspace{.2in}
\begin{tabular}{| c |} \hline
1 \\ \hline
3 \\ \hline \hline
2 \\ \hline
\end{tabular} \hspace{.2in}
\begin{tabular}{| c |} \hline
2 \\ \hline
1 \\ \hline \hline
3 \\ \hline
\end{tabular} \hspace{.2in}
\begin{tabular}{| c |} \hline
2 \\ \hline
3 \\ \hline \hline
1 \\ \hline
\end{tabular} \hspace{.2in}
\begin{tabular}{| c |} \hline
3 \\ \hline
1 \\ \hline \hline
2 \\ \hline
\end{tabular} \hspace{.2in}
\begin{tabular}{| c |} \hline
3 \\ \hline
2 \\ \hline \hline
1 \\ \hline
\end{tabular}
\end{center}
\caption{Tableaux with cutoffs for approval voting.}
\label{figure:  tableaux with cutoffs}
\end{figure}

Although it is perfectly fine for a voter to approve of all or none of the candidates in approval voting, such a preference will have absolutely no impact on the outcome of the election.  For convenience, we will therefore assume that our voters approve of at least one, but not all, of the candidates.  In this setting, our \emph{ranked approval profile} looks like
\[
\mathbf{p} = 
\begin{bmatrix}
\mathbf{p}_1 \\
\mathbf{p}_2 \\
\vdots \\
\mathbf{p}_{n-1}
\end{bmatrix}
\]
where $\mathbf{p}_i$ corresponds to those voters who have approved of exactly $i$ candidates.  

Approval voting is related to positional voting in that we may view it as being made up of several positional voting systems occurring simultaneously.  To explain, consider the weighting vector
\[
\mathbf{a}_i = (\underbrace{1, \dots, 1}_{i}, \underbrace{0, \dots, 0}_{n-i})
\]
and note that the sum-zero portion of the results of the election using approval voting is given by
\[
\mathbf{p}_1 \widehat{\mathbf{a}}_1 + \mathbf{p}_2 \widehat{\mathbf{a}}_2 + \cdots + \mathbf{p}_{n-1} \widehat{\mathbf{a}}_{n-1} = \mathbf{r}_{\textnormal{app}}.
\]
On the other hand, for a positional vote with respect to the weighting vector $\mathbf{w}$, the sum-zero portion of the result is given by 
\[
(\mathbf{p}_1 + \mathbf{p}_2 + \cdots + \mathbf{p}_{n-1}) \widehat{\mathbf{w}} = \mathbf{r}_{\textnormal{pos}}.
\]

In the spirit of extending the ideas found in Theorem~\ref{theorem:  big one} to other settings such as approval voting, the following theorem shows that $\mathbf{r}_{\textnormal{app}}$ and $\mathbf{r}_{\textnormal{pos}}$ need not have anything in common.  

\begin{theorem}\label{theorem:  approval voting}
Let $n \ge 3$, let $\mathbf{r}_{\textnormal{app}}$ and $\mathbf{r}_{\textnormal{pos}}$ be any two sum-zero results vectors in $\mathbb{Q}^n$, and let $\mathbf{w}$ be any nontrivial sum-zero weighting vector in $\mathbb{Q}^n$.  Then there exist infinitely many ranked approval profiles
\[
\mathbf{p} = 
\begin{bmatrix}
\mathbf{p}_1 \\
\mathbf{p}_2 \\
\vdots \\
\mathbf{p}_{n-1}
\end{bmatrix}
\]
such that the approval voting outcome of $\mathbf{p}$ is $\mathbf{r}_{\textnormal{app}}$, and the positional voting outcome with respect to $\mathbf{w}$ is $\mathbf{r}_{\textnormal{pos}}$. 
\end{theorem}

\begin{proof}
It must be the case that $\mathbf{w}$ and $\widehat{\mathbf{a}}_1$ are linearly independent, or $\mathbf{w}$ and $\widehat{\mathbf{a}}_2$ are linearly independent.  Suppose, without loss of generality, that $\mathbf{w}$ and $\widehat{\mathbf{a}}_2$ are linearly independent.  Set $\mathbf{p}_1$ and $\mathbf{p}_2$ such that $\mathbf{p}_1 \widehat{\mathbf{a}}_1 = \mathbf{r}_{\textnormal{app}}$, $\mathbf{p}_2 \widehat{\mathbf{a}}_2 = \mathbf{0}$, and $\mathbf{p}_2 \mathbf{w} = \mathbf{r}_{\textnormal{pos}} - \mathbf{p}_1 \mathbf{w}$, which we know can be done by Theorem~\ref{theorem:  big one}.  Then set $\mathbf{p}_3 = \cdots = \mathbf{p}_{n-1} = \mathbf{0}$.  The resulting ranked approval profile $\mathbf{p}$ has the desired property.  Furthermore, by Theorem~\ref{theorem:  big one}, there are an infinite number of such ranked approval profiles.
\end{proof}


\section{Equivalent Weighting Vectors and Effective Spaces}

At this point, it is helpful to put an equivalence relation on weighting vectors.  The idea is that two weighting vectors should be equivalent if and only if they yield the same ordinal rankings for all profiles.  For convenience, we focus our attention on the fully ranked situation throughout this section.  

Let the all-ones vector in $\mathbb{Q}^n$ be denoted by $\mathbf{1}$.  We say that two weighting vectors $\mathbf{w}$ and $\mathbf{x}$ in $\mathbb{Q}^n$ are equivalent, and write $\mathbf{w} \sim \mathbf{x}$, if and only if there exist $\alpha, \beta \in \mathbb{Q}$ such that $\alpha > 0$ and $\mathbf{x} = \alpha \mathbf{w} + \beta \mathbf{1}$.  This equivalence relation is often used in the literature to simplify calculations and to pinpoint non-cosmetic differences between different positional voting procedures.  

To motivate this equivalence relation, note that, for any positive rational number $\alpha \in \mathbb{Q}$, it makes sense to say that $\mathbf{w}$ is equivalent to $\alpha \mathbf{w}$ since, for every $\mathbf{p} \in M^{(1, \dots, 1)}$, the ordinal ranking given by $T_{\mathbf{w}}(\mathbf{p})$ is the exactly the same as that of $T_{\alpha \mathbf{w}}(\mathbf{p})$.  After all, the entries in $T_{\alpha \mathbf{w}}(\mathbf{p})$ are simply the entries of $T_{\mathbf{w}}(\mathbf{p})$ multiplied by $\alpha > 0$.  

Furthermore, suppose $\alpha, \beta \in \mathbb{Q}$ where $\alpha > 0$.  If $\mathbf{x} =  \alpha \mathbf{w} + \beta \mathbf{1}$, then the ordinal ranking given by $T_{\mathbf{x}}(\mathbf{p})$ is exactly the same as that given by $T_{\mathbf{w}}(\mathbf{p})$.  This is because the addition of $\beta \mathbf{1}$ to $\alpha \mathbf{w}$ changes each candidate's score by exactly the same amount.  

Note that $\mathbf{w} \sim \mathbf{x}$ if and only if there is a positive rational number $\gamma \in \mathbb{Q}$ such that $\widehat{\mathbf{w}} = \gamma \widehat{\mathbf{x}}$, i.e., the sum-zero component of $\mathbf{w}$ is a positive multiple of the sum-zero component of $\mathbf{x}$.  This is helpful to see because, by Theorem~\ref{theorem:  big one}, it means that two weighting vectors $\mathbf{w}$ and $\mathbf{x}$ will always yield the same outcome if and only if $\mathbf{w} \sim \mathbf{x}$.  We therefore have the following theorem:

\begin{theorem}\label{theorem:  ordinal equiv}\textnormal{(Theorem 2.3.1 in~\cite{saari-1994-geometry})}
Let $n \ge 2$, and let $\mathbf{w}$ and $\mathbf{x}$ be weighting vectors in $\mathbb{Q}^n$.  The ordinal rankings of $T_\mathbf{w} (\mathbf{p})$ and $T_\mathbf{x} (\mathbf{p})$ will be the same for all $\mathbf{p} \in M^{(1, \dots, 1)}$ if and only if $\mathbf{w} \sim \mathbf{x}$.  
\end{theorem}

It is helpful to view Theorem~\ref{theorem:  ordinal equiv} in terms of \emph{effective spaces}, which is an approach used extensively and with great success by Saari (see, for example, \cite{saari-2000-mathematical-1}, \cite{saari-2000-mathematical-2}, \cite{saari-2001-chaotic}, and \cite{saari-2002-adopting}).  Recall that the effective space $E(T)$ of a linear transformation $T$ is the orthogonal complement $\ker (T)^\perp$ of the kernel of $T$.  For convenience, if $\mathbf{w}$ is a weighting vector, then we will denote the effective space of $T_\mathbf{w}$ by $E(\mathbf{w})$ (rather than $E(T_\mathbf{w})$).

As a $\mathbb{Q} S_n$-submodule of the profile space $\mathbb{Q} S_n$, the effective space $E(\mathbf{w})$ of any nontrivial sum-zero weighting vector $\mathbf{w}$ is isomorphic to $S^{(n-1,1)}$.  If $\mathbf{w} \in \mathbb{Q}^n$ has a nontrivial projection onto the all-ones vector and $\widehat{\mathbf{w}} \ne \mathbf{0}$, then $E(\mathbf{w}) \cong  S^{(n)} \oplus S^{(n-1,1)}$.  On the other hand, if $\mathbf{w}$ is simply a nonzero multiple of the all-ones vector (so $\widehat{\mathbf{w}} = \mathbf{0}$), then $E(\mathbf{w}) \cong  S^{(n)}$ (and we only get ties).  

\begin{theorem}\label{theorem:  nonequiv empty intersection}
Let $\mathbf{w}$ and $\mathbf{x}$ be nontrivial sum-zero weighting vectors in $\mathbb{Q}^n$.  Then $E(\mathbf{w}) = E(\mathbf{x})$ if and only if $\mathbf{w} \sim \mathbf{x}$ or $\mathbf{w} \sim - \mathbf{x}$.  Furthermore, if $E(\mathbf{w}) \ne E(\mathbf{x})$, then $E(\mathbf{w}) \cap E(\mathbf{x}) = \{ \mathbf{0} \}$.  
\end{theorem}

\begin{proof}
Suppose $E(\mathbf{w}) = E(\mathbf{x})$.  This implies that $\ker(T_{\mathbf{w}}) = \ker(T_{\mathbf{x}}) = \ker(T_{- \mathbf{x}})$.  By Theorem~\ref{theorem:  big one}, if $\mathbf{w} \nsim \mathbf{x}$ and $\mathbf{w} \nsim - \mathbf{x}$, then there exists a profile $\mathbf{p}$ such that $\mathbf{p} \mathbf{w} \ne \mathbf{0}$ and $\mathbf{p} \mathbf{x} = \mathbf{p} (- \mathbf{x}) = \mathbf{0}$.  Thus, if $E(\mathbf{w}) = E(\mathbf{x})$, then $\mathbf{w} \sim \mathbf{x}$ or $\mathbf{w} \sim - \mathbf{x}$.  

On the other hand, if $\mathbf{w} \sim \mathbf{x}$ or $\mathbf{w} \sim - \mathbf{x}$, then $\mathbf{w}$ and $\mathbf{x}$ are linearly dependent (since we are assuming that $\mathbf{w} = \widehat{\mathbf{w}}$ and $\mathbf{x} = \widehat{\mathbf{x}}$).  Thus, $\ker(T_{\mathbf{w}}) = \ker(T_{\mathbf{x}})$, implying that $E(\mathbf{w}) = E(\mathbf{x})$.  

Finally, $E(\mathbf{w}) \cap E(\mathbf{x})$ is a submodule of both $E(\mathbf{w})$ and $E(\mathbf{x})$, and $E(\mathbf{w})$ and $E(\mathbf{x})$ are simple submodules (that are isomorphic to $S^{(n-1,1)}$).  Thus, if $E(\mathbf{w}) \ne E(\mathbf{x})$, then it follows that $E(\mathbf{w}) \cap E(\mathbf{x}) = \{ \mathbf{0} \}$.
\end{proof}

By Theorem~\ref{theorem:  nonequiv empty intersection}, distinct effective spaces for sum-zero weighting vectors intersect only at $\mathbf{0}$.  We can, however, say more.  To explain, we write $\mathbf{w} \perp \mathbf{x}$ if the dot product of $\mathbf{w}$ and $\mathbf{x}$ is zero, i.e., if they are orthogonal.  Furthermore, if $U$ and $W$ are subspaces of a vector space such that every vector in $U$ is orthogonal to every vector in $W$, then we write $U \perp W$.  

Recall that we may view permutations $\sigma \in S_n$ as tableaux in $M^{(1, \dots, 1)}$.  For example, the permutation $\sigma = (1 2 4)(3 5)$ corresponds to the tableau
\begin{center}
\begin{tabular}{| c |} \cline{1-1}
2 \\ \hline 
4 \\ \hline
5 \\ \hline
1 \\ \hline
3 \\ \hline
\end{tabular}
\end{center}
in $M^{(1,1,1,1,1)}$.  In particular, note that the position that candidate $j$ occupies with respect to the permutation $\sigma$ is given by $\sigma^{-1}(j)$, in which case, for a weighting vector $\mathbf{w} = [w_1, \dots, w_n]^t$, candidate $j$ would receive $w_{\sigma^{-1}(j)}$ points.  

\begin{theorem}
If $\mathbf{w}$ and $\mathbf{x}$ are nontrivial sum-zero weighting vectors in $\mathbb{Q}^n$, then $E(\mathbf{w}) \perp E(\mathbf{x})$ if and only if $\mathbf{w} \perp \mathbf{x}$.  
\end{theorem}

\begin{proof}
First, note that $E(\mathbf{w})$ is simply the row space of $T_\mathbf{w}$ when we view $T_\mathbf{w}$ as a matrix with respect to the indicator functions of $M^{(1, \dots, 1)}$.  It follows that if $E(\mathbf{w}) \perp E(\mathbf{x})$, then each row of $T_{\mathbf{w}}$ is orthogonal to each row of $T_{\mathbf{x}}$.  The dot product of the first row of $T_{\mathbf{w}}$ and the first row of $T_{\mathbf{x}}$, however, is a non-zero multiple of the dot product of $\mathbf{w}$ and $\mathbf{x}$, as we show below in~(\ref{equation:  useful result}).  It follows that, if $E(\mathbf{w}) \perp E(\mathbf{x})$, then $\mathbf{w} \perp \mathbf{x}$.   

On the other hand, suppose $\mathbf{w} \perp \mathbf{x}$.  Partition the permutations of the candidates into $n$ sets $X_1, \dots, X_n$ where $X_i$ contains the permutations that have the first candidate, $c_1$, in the $i$th position.  Within each $X_i$, every candidate other than $c_1$ occupies every position other than the $i$th position the same number of times, namely $(n-2)!$ times.  This is because the $i$th position is taken by $c_1$, and by fixing $c_j$, $j \ne 1$, in some position, we are free to place the other candidates in $(n-2)!$ ways.  

The rows of $T_{\mathbf{w}}$ correspond to functions defined on the permutations of the $n$ candidates.  The value that the $j$th row assigns to the permutation $\sigma$ is $w_{\sigma^{-1}(j)}$, which is, of course, the weight given to candidate $j$ based on the permutation $\sigma$.  

Let $r_1(\mathbf{w})$ be the first row of $T_{\mathbf{w}}$, and let $r_j(\mathbf{x})$ be the $j$th row of $T_{\mathbf{x}}$.  These rows may be viewed as elements of $\mathbb{Q}S_n$, where $(r_1(\mathbf{w}))(\sigma) = w_{\sigma^{-1}(1)}$ and $(r_j(\mathbf{x}))(\sigma) = x_{\sigma^{-1}(j)}$.  Taking dot products yields 
\begin{align*}
r_1(\mathbf{w}) \cdot r_j(\mathbf{x}) 
&= \sum_{\sigma} w_{\sigma^{-1}(1)} x_{\sigma^{-1}(j)} \\
&= \sum_{i=1}^n \sum_{\sigma \in X_i} w_{\sigma^{-1}(1)} x_{\sigma^{-1}(j)} \\
&= \sum_{i=1}^n w_i \sum_{\sigma \in X_i} x_{\sigma^{-1}(j)}.
\end{align*}

If $j = 1$, then 
\[
\sum_{\sigma \in X_i} x_{\sigma^{-1}(j)} = \sum_{\sigma \in X_i} x_{\sigma^{-1}(1)}
= (n-1)! x_i.  
\]
It follows that 
\begin{align}\label{equation:  useful result}
r_1(\mathbf{w}) \cdot r_1(\mathbf{x}) 
&= \sum_{i=1}^n w_i ((n-1)! x_i) \\
&= (n-1)! \sum_{i=1}^n w_i x_i \notag \\
&= (n-1)! (\mathbf{w} \cdot \mathbf{x}) \notag \\
&= 0  \notag
\end{align}
since $\mathbf{w} \perp \mathbf{x}$.  

On the other hand, if $j \ne 1$, then 
\[
\sum_{\sigma \in X_i} x_{\sigma^{-1}(j)} 
= (n-2)! \sum_{k \ne i} x_k
= (n-2)! (- x_i).    
\]
which implies that
\begin{align*}
r_1(\mathbf{w}) \cdot r_j(\mathbf{x}) 
&= \sum_{i=1}^n w_i ((n-2)! (- x_i)) \\
&= -(n-2)! \sum_{i=1}^n w_i x_i \\
&= -(n-2)! (\mathbf{w} \cdot \mathbf{x}) \\
&= 0.  
\end{align*}
Thus, if $\mathbf{w} \perp \mathbf{x}$, then we have that $r_1(\mathbf{w}) \perp r_j(\mathbf{x})$ for all $1 \le j \le n$.  To complete the proof, note that the $i$th row of $T_\mathbf{w}$ is the result of acting on $r_1(\mathbf{w})$ with the transposition $\xi = (1\ i)$ that swaps 1 and $i$.  It follows that 
\[
r_i(\mathbf{w}) \cdot r_j(\mathbf{x}) 
= (\xi r_i(\mathbf{w})) \cdot (\xi r_j(\mathbf{x}))
= r_1(\mathbf{w}) \cdot r_{\xi(j)}(\mathbf{x}) = 0.
\]
Thus, the row space of $T_{\mathbf{w}}$ is orthogonal to the row space of $T_{\mathbf{x}}$.  In other words, $E(\mathbf{w}) \perp E(\mathbf{x})$.  
\end{proof}


\section{The Borda Count}

If you are familiar at all with Saari's work, then you know that the Borda count, i.e., the positional voting procedure for $n$ candidates that uses the weighting vector $\mathbf{w} = [n-1, n-2, \dots, 2, 1, 0]^t$, plays a special role when it comes to positional voting.  In this section, we use the algebraic framework we have thus far created to show why this is the case.  In doing so, we also begin to pave the way toward an analogue to the Borda count for partially ranked voting data.  

To motivate our discussion, consider the so-called \emph{Copeland Method} for running an election.  This procedure is based on information concerning head-to-head contests between the candidates.  For each candidate $c_i$, let $w(i)$ and $l(i)$ be the number of head-to-head contests won and lost, respectively, by $c_i$.  The winner under Copeland's Method is the candidate whose difference $w(i) - l(i)$ is largest.  

In our running example of an election with three candidates and eleven voters, candidate $c_3$ defeats both $c_1$ and $c_2$ in head to head contests, and $c_2$ defeats $c_1$.  The scores for candidates $c_1$, $c_2$, and $c_3$ are therefore $0-2=-2$, $1-1=0$, and $2-0=2$, respectively.  Thus, $c_3$ is the winner using the Copeland Method.  Note, by the way, that $c_3$ beat all of the other candidates in head-to-head contests.  When such a candidate exists, she is said to be a \emph{Condorcet winner}.  

What makes the Copeland Method interesting for us is that all of the results can be derived from the image of a map $P:  M^{(1, \dots, 1)} \to M^{(1,1,n-2)}$ which we call the \emph{pairs map}.   The idea behind the pairs map is that it extracts all of the necessary information concerning pairs of candidates (think head-to-head contests).  The defining characteristic of $P$ is that it maps a basis vector $\mathbf{u}$ in $M^{(1, \dots, 1)}$ to the sum of all basis vectors in $M^{(1,1,n-2)}$ whose associated ordered pairs are ranked in the same order as they are ranked in $\mathbf{u}$.  

For example, suppose $n = 4$.  Then the image of the basis vector corresponding to 
\begin{center}
\begin{tabular}{| c |} \hline
3 \\ \hline
1 \\ \hline
4 \\ \hline
2 \\ \hline
\end{tabular} 
\end{center}
is 
\begin{center}
\begin{tabular}{| c c |} \cline{1-1}
3 & \multicolumn{1}{| c}{} \\ \cline{1-1}
1 & \multicolumn{1}{| c}{} \\ \cline{1-2}
2 & 4 \\ \cline{1-2}
\end{tabular} +
\begin{tabular}{| c c |} \cline{1-1}
3 & \multicolumn{1}{| c}{} \\ \cline{1-1}
4 & \multicolumn{1}{| c}{} \\ \cline{1-2}
1 & 2 \\ \cline{1-2}
\end{tabular} +
\begin{tabular}{| c c |} \cline{1-1}
3 & \multicolumn{1}{| c}{} \\ \cline{1-1}
2 & \multicolumn{1}{| c}{} \\ \cline{1-2}
1 & 4 \\ \cline{1-2}
\end{tabular} +
\begin{tabular}{| c c |} \cline{1-1}
1 & \multicolumn{1}{| c}{} \\ \cline{1-1}
4 & \multicolumn{1}{| c}{} \\ \cline{1-2}
2 & 3 \\ \cline{1-2}
\end{tabular} +
\begin{tabular}{| c c |} \cline{1-1}
1 & \multicolumn{1}{| c}{} \\ \cline{1-1}
2 & \multicolumn{1}{| c}{} \\ \cline{1-2}
3 & 4 \\ \cline{1-2}
\end{tabular} +
\begin{tabular}{| c c |} \cline{1-1}
4 & \multicolumn{1}{| c}{} \\ \cline{1-1}
2 & \multicolumn{1}{| c}{} \\ \cline{1-2}
1 & 3 \\ \cline{1-2}
\end{tabular}.
\end{center}

Given the profile $\mathbf{p} \in M^{(1, \dots, 1)}$, the scores for Copeland's Method can all be determined from the image of $\mathbf{p}$ under the pairs map.  You need only consult the coefficients of $P(\mathbf{p})$ to determine the winner of each head-to-head contest.  Moreover, there are several examples of voting procedures that essentially rely solely on pairs data (see, for example, a list of such procedures in Chapter 4 of~\cite{arrow-sen-suzumura-2002-handbook}). 

An interesting question now arises.  What relationship, if any, is there between a map $T_{\mathbf{w}}$ and the pairs map $P$?  To make this question more concrete, let $T:  V \to W$ and $T':  V \to U$ be two linear transformations defined on the same vector space $V$.  We say that $T'$ is \emph{recoverable} from $T$ if there exists a linear transformation $R:  W \to U$ such that $T' = R \circ T$.  It is easy to show that $T'$ is recoverable from $T$ if and only if $\ker(T) \subseteq \ker(T')$ which, in turn, occurs if and only if $E(T') \subseteq E(T)$.  

This leads us to a much more specific form of the question above.  For what weighting vectors $\mathbf{w} \in \mathbb{Q}^n$ is $T_\mathbf{w}$ recoverable from $P$?  To answer this question, we will focus on the effective spaces of our positional voting procedures and pairs map.  Once again, the  representation theory of the symmetric group will play an important role.

The first thing we want to do is to note that the pairs map $P:  M^{(1, \dots, 1)} \to M^{(1,1,n-2)}$ is a $\mathbb{Q}S_n$-module homomorphism.  Thus, we may make use of Schur's Lemma.  Next, we turn our attention to the effective space of $P$.  The codomain $M^{(1,1,n-2)}$ of the pairs map $P$ has the following decomposition into simple submodules: 
\[
M^{(1,1,n-2)} \cong S^{(n)} \oplus S^{(n-1,1)} \oplus S^{(n-1,1)} \oplus S^{(n-2,2)} \oplus S^{(n-2,1,1)}.
\]
Furthermore, it can be shown (using, for example, a dimension argument) that the image, and therefore the effective space, of $P$ is isomorphic to $S^{(n)} \oplus S^{(n-1,1)} \oplus S^{(n-2,1,1)}$.  Since there is only one copy of $S^{(n-1,1)}$ in this decomposition, it follows by Theorem~\ref{theorem:  nonequiv empty intersection} that there are at most two nontrivial equivalence classes of weighting vectors whose effective spaces are contained in the effective space of $P$.  As the following theorem (which is essentially implied by Theorem 3.2.1 in~\cite{saari-1994-geometry}) shows, there are such equivalence classes.  They are the equivalence classes that contain the Borda count and its negative.  

\begin{theorem}\label{theorem:  borda full}
Let $n \ge 2$, and let $\mathbf{w} \in \mathbb{Q}^n$ be a nontrivial weighting vector (i.e., $\mathbf{w} \nsim \mathbf{1}$).  The map $T_{\mathbf{w}}$ is recoverable from the pairs map $P$ if and only if $\mathbf{w}$ or $- \mathbf{w}$ is equivalent to the Borda count.  
\end{theorem}

\begin{proof}
By the above discussion, it is enough to show that if $\mathbf{w}$ is the Borda count weighting vector, i.e., $\mathbf{w} = [n-1, n-2, \dots, 2, 1, 0]^t \in \mathbb{Q}^n$, then $T_\mathbf{w}$ is recoverable from the $P$.  This, however, is trivial.  In fact, the results vector one obtains by using the Borda count can be (and often is) viewed as the sum of the points awarded to a candidate from all of her head-to-head victories, and these points are encoded (blatantly) in the image of the pairs map $P$.   
\end{proof}

One of the nice properties that the Borda count enjoys is that, if there is a Condorcet winner, she is never ranked last by the Borda count (see, for example, Corollary 5 in~\cite{saari-2000-mathematical-1}).  By Theorem~\ref{theorem:  big one}, any weighting vector that is not equivalent to $\mathbf{w} = [n-1, n-2, \dots, 2, 1, 0]^t$ \emph{does not} enjoy this property.  In fact, in the class of positional voting procedures for fully ranked profiles, the Borda count maximizes the probability that a Condorcet winner is actually ranked first \cite{vannewenhizen-1992-borda}.  For more on the relationship between the Borda count and Copeland's method, see \cite{merlin-saari-1997-copeland-2} and \cite{saari-merlin-1996-copeland-1}.

The Borda count also has what is called \emph{reversal symmetry}.  In other words, under the Borda count, if all of the voters were to completely reverse their ballots so that their first choice is now their last, their second choice is now second to last, and so on, then the resulting ordinal ranking would be the complete reversal of the original result.  When $n = 3$, the Borda count is the unique weighting vector (up to equivalence) with this property, but when $n \ge 4$, there are others.  For example, $\mathbf{w} = [6, 5, 1, 0]^t$ has this property.  This is easy to see, however, once you recognize that $[6, 5, 1, 0]^t \sim [3, 2, -2, -3]^t$.


\section{Analogues to the Borda Count}

Recall that if voters are returning fully ranked ballots, then the Borda count and its negative are the unique (up to equivalence) nontrivial positional voting procedures that are recoverable from the pairs map.  What if, however, the voters do not return fully ranked ballots?  What if it has been decided that it is infeasible to ask voters to rank \emph{all} of the candidates?  

In this section, we turn our attention to the ``rank-only-your-top-$k$" situation in which $\lambda = (1, \dots, 1, n-k) = (1^k, n-k)$.  By generalizing only slightly the pairs map $P:  M^{(1, \dots, 1)} \to M^{(1,1,n-2)}$, we are able to generalize Theorem~\ref{theorem:  borda full} to the ``rank-only-your-top-$k$" situation.  Interestingly, and in contrast to the fully ranked case, we show that there is more than one ``Borda-like" equivalence class of weighting vectors.  

We generalize the pairs map as follows.  Let $0 \le \tau \le 1$, and define $P^k_\tau:  M^{(1^k, n-k)} \to M^{(1,1, n-2)}$ as we did for the pairs map $P$ in the full ranking case, except now, if two candidates $c_i$ and $c_j$ are tied for last place, then we assign both of the ordered pairs $(c_i, c_j)$ and $(c_j, c_i)$ the value $\tau$ (think ``points for tying").  By letting $\tau$ be a parameter, we are able to consider simultaneously an infinite number of analogues of the pairs map $P$.  

For example, suppose $n = 4$ and $\lambda = (1,1,2)$.  Then the image of the basis vector corresponding to 
\begin{center}
\begin{tabular}{| c c |} \cline{1-1}
2 & \multicolumn{1}{| c}{} \\ \cline{1-1}
4 & \multicolumn{1}{| c}{} \\ \cline{1-2}
1 & 3 \\ \cline{1-2}
\end{tabular} 
\end{center}
is 
\begin{center}
\begin{tabular}{| c c |} \cline{1-1}
2 & \multicolumn{1}{| c}{} \\ \cline{1-1}
4 & \multicolumn{1}{| c}{} \\ \cline{1-2}
1 & 3 \\ \cline{1-2}
\end{tabular} +
\begin{tabular}{| c c |} \cline{1-1}
2 & \multicolumn{1}{| c}{} \\ \cline{1-1}
1 & \multicolumn{1}{| c}{} \\ \cline{1-2}
3 & 4 \\ \cline{1-2}
\end{tabular} +
\begin{tabular}{| c c |} \cline{1-1}
2 & \multicolumn{1}{| c}{} \\ \cline{1-1}
3 & \multicolumn{1}{| c}{} \\ \cline{1-2}
1 & 4 \\ \cline{1-2}
\end{tabular} +
\begin{tabular}{| c c |} \cline{1-1}
4 & \multicolumn{1}{| c}{} \\ \cline{1-1}
1 & \multicolumn{1}{| c}{} \\ \cline{1-2}
2 & 3 \\ \cline{1-2}
\end{tabular} +
\begin{tabular}{| c c |} \cline{1-1}
4 & \multicolumn{1}{| c}{} \\ \cline{1-1}
3 & \multicolumn{1}{| c}{} \\ \cline{1-2}
1 & 2 \\ \cline{1-2}
\end{tabular} + $\tau$ \
\begin{tabular}{| c c |} \cline{1-1}
1 & \multicolumn{1}{| c}{} \\ \cline{1-1}
3 & \multicolumn{1}{| c}{} \\ \cline{1-2}
2 & 4 \\ \cline{1-2}
\end{tabular} + $\tau$ \
\begin{tabular}{| c c |} \cline{1-1}
3 & \multicolumn{1}{| c}{} \\ \cline{1-1}
1 & \multicolumn{1}{| c}{} \\ \cline{1-2}
2 & 4 \\ \cline{1-2}
\end{tabular}.
\end{center}

We now have the following question.  For which partial weighting vectors $\mathbf{w} = [w_1, \dots, w_{k+1}]^t$ is $T_\mathbf{w}$ recoverable from $P^k_\tau$?  To answer this question, define $\mathbf{b} = [b_1, \dots, b_{k+1}]^t$ to be the partial weighting vector corresponding to $\lambda = (1^k, n-k)$ where $b_i = n-i$ for $1 \le i \le k$ and 
\[
b_{k+1} = \frac{1}{2} (n-k-1).
\]
This is the partial weighting vector one would get by ``averaging the Borda count with respect to $\lambda$."   In other words, we essentially use the Borda count for the top $k$ candidates, but we assign the average of the last $n-k$ Borda count points to each of the last $n-k$ candidates.  This average is 
\[
\frac{1}{n-k} (0 + 1 + 2 + \cdots +(n-k-1)) = \frac{(n-k-1)(n-k)}{2( n-k)} = \frac{1}{2} (n-k-1) = b_{k+1}.
\]

Similarly, we define $\mathbf{b}^\tau = [b^\tau_1, \dots, b^\tau_{k+1}]^t$ in exactly the same way, except that we set
\[
b^\tau_{k+1} = \tau(n-k-1).
\]
In other words, $b^\tau_i = b_i$ for $1 \le i \le k$, but $b^\tau_{k+1} = 2 \tau b_{k+1}$.  Thus, if $\tau = 1/2$, then $\mathbf{b} = \mathbf{b}^\tau$.  

Our first goal is to show that both $T_\mathbf{b}$ and $T_{\mathbf{b}^\tau}$ are recoverable from $P_\tau^k$.  That is, we want to show that there exist linear transformations $\varphi_{\mathbf{b}}$ and $\psi_{{\mathbf{b}^\tau}}$ such that $T_{\mathbf{b}} = \varphi_{\mathbf{b}} \circ P^k_\tau$ and $T_{{\mathbf{b}^\tau}} = \psi_{{\mathbf{b}^\tau}} \circ P^k_\tau$.  With that in mind, define
\[
\psi_{{\mathbf{b}^\tau}}:  M^{(1,1,n-2)} \to M^{(1,n-1)}
\]
by setting
\[
[\psi_{{\mathbf{b}^\tau}}(v)]_i = \sum_{j \ne i} v^{ij}.
\]
In other words, the result corresponding to the $i$th candidate is determined by summing all of the entries corresponding to the ordered pairs $(i,j)$ where $j \ne i$, i.e.,  all of the pairs in which candidate $c_i$ is beating some other candidate.  

The map $\psi_{{\mathbf{b}^\tau}}$ is easily seen to be a $\mathbb{Q}S_n$-module homomorphism.  If we let $u_e$ be the indicator function corresponding to the tabloid that contains the identiy permutation $e \in S_n$, then we may easily check that 
\[
[(\psi_{{\mathbf{b}^\tau}} \circ P^k_\tau) (u_e)]_i = b^\tau_i.
\]
Together with the fact that $\psi_{{\mathbf{b}^\tau}}$ and $P^k_\tau$ are $\mathbb{Q}S_n$-module homomorphisms, this implies that $T_{{\mathbf{b}^\tau}} = \psi_{{\mathbf{b}^\tau}} \circ P^k_\tau$.  Thus, $T_{{\mathbf{b}^\tau}}$ is recoverable from $P^k_\tau$.  

Similarly, we may construct a linear transformation $\varphi_{\mathbf{b}}$ such that $T_{\mathbf{b}} = \varphi_{\mathbf{b}} \circ P^k_\tau$.  First, note that for any indicator function $u$, if we sum the entries of $P^k_\tau (u)$, we always get the same value $E$ where
\[
E = (n-1) + (n-2) + \cdots + (n-k) + 2t \binom{n-k}{2}.  
\]
With that in mind, we define the $\mathbb{Q}S_n$-module homomorphism
\[
\varphi_{\mathbf{b}}:  M^{(1,1,n-2)} \to M^{(1,n-1)}
\]
by setting
\[
[\varphi_{\mathbf{b}}(v)]_i 
= \frac{1}{2} \left( \sum_{j \ne i} (v^{ij} - v^{ji}) + \frac{n-1}{E} \sum_{k,l} v^{kl} \right) .
\]
Again, we may check (perhaps with a bit more work this time) that 
\[
[(\varphi_{\mathbf{b}} \circ P^k_\tau) (u_e)]_i = b_i.
\]
Thus, $T_{\mathbf{b}} = \varphi_{\mathbf{b}} \circ P^k_\tau$, implying that $T_{\mathbf{b}}$ is also recoverable from $P^k_\tau$.  

Our next goal is to show that $\overline{\mathbf{b}}$ and $\overline{{\mathbf{b}^\tau}}$ are equivalent if and only if $\tau = 1/2$.  This is straightforward.  If they were equivalent, then the ratios of the differences between successive entries would have to be the same.  In particular, it would be the case that 
\[
\frac{b_{k-1} - b_k}{b_k - b_{k+1}} = \frac{b^\tau_{k-1} - b^\tau_k}{b^\tau_k - b^\tau_{k+1}}.
\]
This, however, is true if and only if $\tau = 1/2$, since the above equation reduces to 
\[
\frac{1}{(n-k) - \frac{1}{2}(n-k-1)} = \frac{1}{(n-k) - \tau(n-k-1)}
\]
and solving for $\tau$ shows that $\tau = 1/2$.  We therefore have the following proposition:

\begin{proposition}
The weighting vectors $\overline{\mathbf{b}}$ and $\overline{\mathbf{b}^\tau}$ are equivalent if and only if  $\tau = 1/2$ (in which case $\overline{\mathbf{b}} = \overline{\mathbf{b}^\tau}$).  
\end{proposition}

Finally, we turn our attention to characterizing those weighting vectors that are recoverable from $P^k_\tau$.  We begin with a proposition.  

\begin{proposition}\label{proposition:  one copy}
If $\tau = 1/2$, then the image of $P^k_\tau:  M^{(1^k, n-k)} \to M^{(1,1,n-2)}$ contains exactly one copy of the simple module $S^{(n-1,1)}$.  
\end{proposition}

\begin{proof}
Let $I$ be the image of $P^k_\tau$.  Since $I$ is a submodule of $M^{(1,1,n-2)}$, and the simple module $S^{(n-1,1)}$ appears exactly twice in any decomposition of $M^{(1,1,n-2)}$ into simple modules, we know that any decomposition of $I$ into simple modules can contain at most two copies of $S^{(n-1,1)}$.  Furthermore, since $\mathbf{b}$ is a nontrivial weighting vector and $T_{\mathbf{b}} = \varphi_{\mathbf{b}} \circ P^k_\tau$, we know that $I$ must contain at least one copy of $S^{(n-1,1)}$.  

Let $J$ be the direct sum of the two simple modules that are isomorphic to $S^{(n-1,1)}$ in a decomposition of $M^{(1,1,n-2)}$ into simple modules.  In other words, $J$ is the so-called \emph{isotypic subspace} of $M^{(1,1,n-2)}$ corresponding to $S^{(n-1,1)}$.  It turns out that any simple submodule of $M^{(1,1,n-2)}$ that is isomorphic to $S^{(n-1,1)}$ is necessarily a submodule of $J$.  We will make use of this fact shortly.  

Assume, for the sake of contradiction, that $I$ contains two copies of $S^{(n-1,1)}$ when we decompose it into simple modules.  In other words, assume that $J \subseteq I$.  Since $\tau = 1/2$, we have that $T_{\mathbf{b}} =T _{{\mathbf{b}^\tau}}$.  This implies that $\varphi_{\mathbf{b}} \circ P^k_\tau = \psi_{{\mathbf{b}^\tau}} \circ P^k_\tau$, and therefore that
\[
(\varphi_{\mathbf{b}}  - \psi_{{\mathbf{b}^\tau}}) \circ P^k_\tau
\]
is the zero linear transformation.  It follows that any vector in $I$ must be in $\ker (\varphi_{\mathbf{b}}  - \psi_{{\mathbf{b}^\tau}})$.  In other words, $I \subseteq \ker (\varphi_{\mathbf{b}}  - \psi_{{\mathbf{b}^\tau}})$.

Our assumption is that $J \subseteq I$.  Given the above, we may contradict this assumption by finding a vector in $J$ that is not in $\ker (\varphi_{\mathbf{b}}  - \psi_{{\mathbf{b}^\tau}})$.  With this in mind, for each $1 \le i \le n$, let $\mathbf{v}_i \in M^{(1,1,n-2)}$ be such that the coefficient corresponding to an ordered pair that contains $i$ is $(n-2)/2$, and is $-1$ otherwise.  

The sum of the entries of $\mathbf{v}_i$ is zero.  Furthermore, these vectors generate a submodule that is isomorphic to $S^{(n-1,1)}$.  Thus each $\mathbf{v}_i$ is in the $S^{(n-1,1)}$ isotypic space $J$ of $M^{(1,1,n-2)}$, implying that each $\mathbf{v}_i$ is in $I$.  One may easily verify, however, that $\mathbf{v}_i \in \ker \varphi_{\mathbf{b}}$, but that $\mathbf{v}_i \notin \ker \psi_{{\mathbf{b}^\tau}}$.  Since this is a contradiction, it must be the case that $I$ contains exactly one copy of $S^{(n-1,1)}$.
\end{proof}

The following theorem characterizes those weighting vectors that are recoverable from the map $P^k_\tau$.  More specifically, it says that, with respect to $P^k_\tau$, weighting vectors related to $\mathbf{b}$ and ${\mathbf{b}^\tau}$ form the analogues of the Borda count weighting vector $[n-1, n-2, \dots, 1, 0]^t$ when it comes to the ``rank-only-your-top-$k$" situation.

\begin{theorem}\label{theorem:  borda analogue}
Let $\mathbf{w}$ be a partial weighting vector with respect to $\lambda = (1^k, n-k)$ where $1 \le k \le n-2$.  The positional map $T_{\mathbf{w}}$ is recoverable from the map $P^k_\tau$ if and only if $\widehat{\overline{\mathbf{w}}}$ is a linear combination of $\widehat{\overline{\mathbf{b}}}$ and $\widehat{\overline{{\mathbf{b}^\tau}}}$.  
\end{theorem}

\begin{proof}
If $\widehat{\overline{\mathbf{w}}}$ is a linear combination of $\widehat{\overline{\mathbf{b}}}$ and $\widehat{\overline{{\mathbf{b}^\tau}}}$, then $T_{\mathbf{w}}$ is clearly recoverable from $P^k_\tau$.  On the other hand, suppose that $T_{\mathbf{w}}$ is recoverable from the pairwise map $P^k_\tau$.  If $\mathbf{b} \ne {\mathbf{b}^\tau}$, then $\{ \widehat{\overline{\mathbf{b}}}, \widehat{\overline{{\mathbf{b}^\tau}}} \}$ is a basis for the recoverable sum-zero weighting vectors since the image of $P^k_\tau$ contains at most two copies of the simple module $S^{(n-1,1)}$.  Thus $\widehat{\overline{\mathbf{w}}}$ is a linear combination of $\widehat{\overline{\mathbf{b}}}$ and $\widehat{\overline{{\mathbf{b}^\tau}}}$.  If $\mathbf{b} = {\mathbf{b}^\tau}$, however, then by Proposition~\ref{proposition:  one copy}, $\{ \widehat{\overline{\mathbf{b}}} \}$ is a basis for the recoverable sum-zero weighting vectors.  In either case, the theorem follows.  
\end{proof}

Recall that one of the nice properties of the Borda count is that, if there is a Condorcet winner, then she is never ranked last by the Borda count.  Condorcet winners make sense in the ``rank-only-your-top-$k$" situation as well since the pairs map $P:  M^{(1, \dots, 1)} \to M^{(1,1,n-2)}$ and $P^k_\tau:  M^{(1^k, n-k)} \to M^{(1,1,n-2)}$ have the same codomain.  It turns out that, if a profile $\mathbf{p} \in M^{(1^k, n-k)}$ has a Condorcet winner with respect to $P^k_\tau$, then it has the same Condorcet winner with respect to all maps $P^k_{\tau'}$ where $0 \le \tau' \le 1$.  In other words, the existence of a Condorcet winner does not depend on $\tau$.  

As in the fully ranked case, it also turns out that if a Condorcet winner exists in the ``rank-only-your-top-$k$" situation, then she will never be ranked last under the positional map $T_\mathbf{b}$.  More importantly, given what we have seen so far, it should hardly come as a surprise that this statement is not true for any other weighting vector $\mathbf{w}$ such that $\overline{\mathbf{b}}$ is not equivalent to $\overline{\mathbf{w}}$.  Therefore, if the notion of a Condorcet winner is important to you, then it would certainly be reasonable to say that $\mathbf{b}$ is the unique (up to equivalence) analogue of the usual Borda count.


\section{Acknowledgments}

We gratefully acknowledge the Reed Institute for Applied Statistics at Claremont McKenna College for funding each of the first three authors while part of this work was done.  Special thanks also to Donald Saari and Francis Su for helpful conversations and suggestions.


\end{document}